%somos.tex:
%%a Plain TeX file by Doron Zeilberger (7 pages)

%begin macros

\baselineskip=14pt
\parskip=10pt
\def\halmos{\hbox{\vrule height0.15cm width0.01cm\vbox{\hrule height
  0.01cm width0.2cm \vskip0.15cm \hrule height 0.01cm width0.2cm}\vrule
  height0.15cm width 0.01cm}}
\font\eightrm=cmr8 
\font\eighttt=cmtt8
\magnification=\magstephalf

\def\1{{\overline{1}}}
\def\2{{\overline{2}}}
\parindent=0pt
\overfullrule=0in
\def\Tilde{\char126\relax}
\def\frac#1#2{{#1 \over #2}}
%\headline={\rm  \ifodd\pageno  \RightHead  \else  \LeftHead  \fi}
%\def\RightHead{\centerline{
%Title
%}}
%\def\LeftHead{ \centerline{Doron Zeilberger}}
%end macros
\bf
\centerline
{
How To Generate As Many Somos-Like Miracles as You Wish
}
\rm
\bigskip
\centerline{ {\it Shalosh B. EKHAD and
Doron 
ZEILBERGER}\footnote{$^1$}
{\eightrm  \raggedright
Department of Mathematics, Rutgers University (New Brunswick),
Hill Center-Busch Campus, 110 Frelinghuysen Rd., Piscataway,
NJ 08854-8019, USA.
%\break
{\eighttt zeilberg  at math dot rutgers dot edu} ,
\hfill \break
{\eighttt http://www.math.rutgers.edu/\~{}zeilberg/} .
March 21, 2013.
Accompanied by Maple package \hfill \break {\eighttt NesSomos}
downloadable from Zeilberger's website.
Supported in part by the NSF.
}
}

{\bf Somos Sequences}

In the late 1980s, Michael Somos came up with an amazing 
set of conjectures, featuring non-linear recurrences that always yield integers.
The most famous one was Somos-4: [S]{\tt http://oeis.org/A006720}. 
Define a sequence $\{ a(n) \}_{n=1}^{\infty}$ by $a(1)=1 \, , \, a(2)=1 \, , \, a(3)=1 \, , \, a(4)=1$, and for $n \geq 5$,
$$
a(n)=\frac{a(n-1)a(n-3)+a(n-2)^2}{a(n-4)}  \quad . 
$$
Prove that $a(n)$ is  an integer for every $n \geq 1$. This
is indeed amazing, since, starting at $n=9$, you divide by an integer larger than $1$, and {\it a-priori}
you are only guaranteed to get {\it rational numbers}, but surprise-surprise, you {\it always}
(seem) to get integers.

This took quite some effort
to prove {\it rigorously}. Of course, confirming it for $n \leq 300$, say, already constitutes
a (very convincing!) {\it empirical} proof, but mathematicians do not care about
{\it truth} per-se, they only care about playing their (sometimes fun but often very dull) {\it artificial} {\bf game}
called  {\it proving} (``rigorously'', using {\it logical} deduction).

And indeed, such a mathematical (elegant!) proof was given by Janice Malouf[M].
Other people gave other proofs, see the references given at
[S] {\tt http://oeis.org/A006720}. 

Of course, Michael Somos is a {\bf honest} person. He first encountered this sequence in
some context, and then  made his conjecture. But suppose that you are {\bf dishonest},
and want to amaze your friends (or challenge your foes) with such {\it Somos Phenomena},
coming up with non-linear recurrences
where you know {\it beforehand}, that you are {\it guaranteed}  to get integers, and even have a fully rigorous proof of that fact.

Then you must follow the advise of Carl Gustav Jacob Jacobi who taught us that in order to make life simpler,

``{\it man muss immer umkehren}'' \quad,

in other words, work backwards! Instead of trying to play darts, by placing the target on the wall, and then
aiming from a distance of twenty meters, and most probably missing it, one {\bf first} throws the dart, and {\it then} draws
the bull's eye!

Consider the linear recurrence sequence defined by
$$
f(1)=1 \, , \, f(2)=3 \, , \, f(n)=3f(n-1)-f(n-2) \quad (n \geq 3) \quad .
$$
Obviously (by induction on $n$) they are {\bf all} integers.

Consider, on the other hand, the following sequence
$$
b(1)=1 \, , \,  b(2)=3 \, , \, b(n)=\frac{b(n-1)^2-1}{b(n-2)} \quad (n \geq 3) \quad ,
$$
now it is not at all immediately obvious that the $b(n)$ are always integers.

Suppose that one asked the question ``prove that the $b(n)$ are always integers''.
The {\bf best}  way to prove it would be to generate the first few terms, go to Neil Sloane's OEIS[S],
and find out that $b(n)=A001519(n)$. In other words it is the same as the
above sequence $f(n)$, at least for the first $30$ terms.

So far it is only a ``conjecture'', but if you know about the
C-finite ansatz ([K][KP][Z1][Z2]), then knowing that $f(n)=b(n)$ for $1 \leq n \leq 4$
implies that it is true for {\it all} $n \geq 1$, thereby yielding an `empirical',
{\bf yet fully rigorous}, proof!

It is easy to see that  C-finite sequences, i.e. sequences that satisfy {\bf linear} recurrence equations with
{\it constant coefficients}, also satisfy (many!) {\bf non-linear} recurrences. But more surprisingly,
the same is true for any {\bf subsequence} where the {\bf places} are either {\bf polynomials} (e.g. $n^2$) or
{\bf exponential} (e.g. $2^n$), or even {\bf exponential polynomial} (e.g. $n^2 2^n+n^4 3^n$). For example, if $F(n)$ are the Fibonacci
numbers defined by $F(0)=0,F(1)=1, F(n)=F(n-1)+F(n-2) (n \geq 2)$, then
the sequence
$$
a(n)\, := \, F(n^2 2^n+ n^3 5^n+1) \quad
$$
satisfies some (very complicated!) non-linear recurrence equation, i.e. there exists a positive integer
$r$ (the order) and a polynomial $P(x_0, \dots , x_r)$ such that
$$
P(a(n), \dots, a(n+r))  = 0 \quad (for \,\,\, all \,\,\ n \geq 1).
$$
Furthermore, by possibly increasing $r$, one can demand that the degree in $x_r$ is $1$, so
one can express $a(n+r)$ as a {\it rational function} of $a(n), \dots, a(n+r-1)$ and get
a genuine non-linear recurrence that {\it very surprisingly} 
(if you don't know how it was found!) yields integers.

Of course the same is true for {\it any} sequence that satisfies a linear recurrence equation with
constant coefficients.

So let's state the main theorem in its  intimidating full generality.

{\bf Main Theorem}: Let $f(n)$ be a sequence satisfying a linear recurrence equation with
constant coefficients
$$
f(n)=\sum_{i=1}^d c_i f(n-i) \quad,
$$
for some constants $c_1, \dots , c_d$, and with some given initial conditions
$f(1), \dots, f(d)$.
Let $p(n)$ be of the form, for some specific integer $m$,
$$
p(n):=POL(n, 2^n,3^n, \dots , m^n) \quad,
$$
where $POL$ is some polynomial of $m$ variables.
Then there exists an integer $r$ and a polynomial $P(x_0, \dots, x_r)$, of degree $1$ in $x_r$ such that
the subsequence \hfill\break $a(n):=f(p(n))$ satisfies the non-linear recurrence
$$
P(a(n), \dots, a(n+r))=0 \quad (n \geq 1 ) \quad  .
$$

The formal proof of this theorem would be very boring, but its informal proof, at least {\it generically}, is not
too bad. Using Binet's formula, one can express $f(n)$  ``explicitly'' as a linear combination of
$\alpha^n$ (and possibly $n \alpha^n, n^2 \alpha^n$ etc.  in case of multiple roots) where
the $\alpha$'s are roots of the equation
$$
1=\sum_{i=1}^d c_i \alpha^{-i}  \quad.
$$
By using the initial conditions, one can find the coefficients in that linear combination of the
$\alpha^n$'s, that are also certain expressions in the
$\alpha$'s. Then, since $\alpha^{(n+1)^2}=\alpha^{n^2} (\alpha^{n})^2 \alpha$, 
$\alpha^{2^{n+1}}=(\alpha^{2^n})^2$ etc. one can introduce {\it auxiliary variables},
$\{ \alpha^n, \alpha^{n^2},\alpha^{n^3}, \dots ,\alpha^{2^n}, \alpha^{3^n} , \dots \}$,
(finitely many!), and  recall that the $\alpha$'s all satisfy the algebraic equation above
(and enjoy nice relations amongst themselves, e.g. that their sum is $c_1$ and their product
is $(-1)^{d-1}c_d$).

Now one can express $f(p(n))$ in terms of these auxiliary variables, and
by  {\it elimination} (in principle using  the Buchberger algorithm), one can get,
for some $r$, a pure polynomial relation, of the desired form
$P(a(n), \dots, a(n+r))=0$, and even demand that the degree in $a(n+r)$ is one, so that one
can express $a(n+r)$ as a rational function of $a(n), \dots , a(n+r-1)$.

In practice, it is better to use {\it undetermined coefficients}. Fixing $r$ and $d$, write
a {\it generic} polynomial $P(x_0, \dots , x_r)$ of total degree $d$ (and degree $1$ in $x_r$), plug-in the
expressions in terms of the auxiliary variables, get a huge polynomial in these variables,
set up all the coefficients (w.r.t. the auxiliary variables)
of the numerator to $0$, and solve the large system of linear equations where the unknowns are the
``undetermined coefficients'', thereby hopefully making them determined. If the only solution
is the trivial zero solution, don't give up! Just raise $d$ and/or $r$.

{\bf Now} if we take the coefficients of the recurrence $c_1, \dots, c_d$ to be {\it integers},
and the initial conditions $f(1), \dots, f(d)$ to be integers as well, then the sequence
$f(n)$ obviously consists of only integers. If, in addition, the polynomial $POL$ has all 
positive integer coefficients, then we are guaranteed, {\it a priori}, that the terms of the subsequence
$a(n):=f(p(n))$ are all integers!

For example, if our sequence is $f(n)=F_{2n}$, i.e. $f(1)=1 \, , \, f(2)=3 \, , \,$ \hfill\break $f(n)=3f(n-1)-f(n-2) \quad (n \geq 3)$,
considered above, then the Binet formula is
$$
f(n)=\frac{\alpha^n- \alpha^{-n}}{\alpha- \alpha^{-1}} \quad,
$$
where $\alpha$ is a root of the quadratic equation $\alpha^2-3\alpha+1=0$. Taking $p(n)=n^2$ we get
$$
a(n)=\frac{\alpha^{n^2}- \alpha^{-n^2}}{\alpha- \alpha^{-1}} \quad .
$$
Now
$$
a(n+1)=\frac{\alpha^{(n+1)^2}- \alpha^{-(n+1)^2}}{\alpha- \alpha^{-1}}=
\frac{\alpha^{n^2} (\alpha^n)^2 \alpha - \alpha^{-n^2} (\alpha^{-n})^2 \alpha^{-1} }{\alpha- \alpha^{-1}} \quad
$$
$$
a(n+2)=\frac{\alpha^{(n+2)^2}- \alpha^{-(n+2)^2}}{\alpha- \alpha^{-1}}=
\frac{\alpha^{n^2} (\alpha^n)^4 \alpha^4 - \alpha^{-n^2} (\alpha^{-n})^4 \alpha^{-4} }{\alpha- \alpha^{-1}} \quad .
$$
So you have expressions for $a(n)$,$a(n+1),a(n+2)$, as rational functions in the auxiliary variables $A:=\alpha^n$, $B:=\alpha^{n^2}$
(and $\alpha$, but $\alpha$ satisfies $\alpha^2-3 \alpha +1=0$). So generically we should have a
second-order polynomial relation $P(a(n),a(n+1),a(n+2))=0$. But it is not of  degree one
in $a(n+2)$. If we allow higher-order recurrences, 
it turns out that we have a fairly nice non-linear recurrence of order $5$. Here is one such example.

{\bf Prop.} Let $a(n)$ be defined by
$$
a(1) = 1 \quad , \quad a(2) = 21 \quad , \quad a(3) = 2584 \quad , \quad a(4) = 2178309 \quad , \quad a(5) = 12586269025
$$
and for $n \geq 6$, by the recurrence 
$$
a(n) = ( \, 2303 \, a(n-4)a(n-3)a(n-1)+2255 \, a(n-3)^2a(n-2) \, + \, 329a(n-4)a(n-1)^2
$$
$$
- \, 15792 \, a(n-4)a(n-2)^2 \, + \, 329 \, a(n-4)a(n-3)^2-2303 \, a(n-4)^2a(n-2)+441 \, a(n-2)
$$
$$
- \, a(n-2)^3-2961 \, a(n-4) \, - \, a(n-5)a(n-2)a(n-1) \, + \, 329 \, a(n-5)a(n-3)a(n-2) \, )/( \, 48 \, a(n-4)a(n-2) \, ) \quad,
$$
then $a(n)$ are all integers!

For many more such propositions, see: \hfill\break
{\tt http://www.math.rutgers.edu/\~{}zeilberg/tokhniot/oNesSomos3} $\,\,$ .

The much sparser subsequence,
$a(n):=f(2^n)$ satisfies a much simpler (second-order!) non-linear recurrence:
$$
a(n)=
{\frac {a \left( n-1 \right)  \left( 4+5\,  a \left( n-1 \right) ^{2} \right) }{2+5\,  a \left( n-2 \right)^{2}}}
$$
subject to the initial conditions
$$
a(1)=3 \quad, \quad a(2)= 21 \quad .
$$

{\bf The Maple package NesSomos}

Everything is implemented in the Maple package {\tt NesSomos} (written by DZ) available directly from

{\tt http://www.math.rutgers.edu/\~{}zeilberg/tokhniot/NesSomos} $\,$,

and that is linked to from the ``front'' of this article

{\tt http://www.math.rutgers.edu/\~{}zeilberg/mamarim/mamarimhtml/somos.html }, where
one can also find three long webbbooks generated by SBE.

To get a list of the main procedures, type {\tt ezra();}, and for help with a specific procedure,
type {\tt ezra(ProcedureName);}. Let's just list the more important procedures.

$\bullet$  {\tt FindREL(LPols,Vars,x,d)} is a very general procedure that inputs a list of rational functions
{\tt LPols} in the list of variables {\tt Vars}, a symbol, {\tt x}, and a positive integer {\tt d}, and outputs a 
polynomial of degree $\leq d$, let's call it $P(x[1], \dots, x[r])$ such that 
$$
P(LPols[1], \dots, LPols[r])= 0 \quad.
$$
For example, 

$$
FindREL( [p^2-q^2,2pq,p^2+q^2],[p,q],x,2); 
$$
yields 
$$
{x_{{3}}}^{2}-{x_{{2}}}^{2}-{x_{{1}}}^{2}=0 \quad .
$$
If nothing is found, then it returns $FAIL$ and one has to increase $d$. Of course 
in order to {\it guarantee} that there is such a good $d$,
we must insist that the length of {\tt LPols} is larger than the length of {\tt Vars}.

$\bullet$  {\tt FindInREL(LPols,Vars,x,d)} is exactly as above, but, in addition,
the output should have
degree one in $x[r]$ (where $r$ is the length of {\tt LPols})). This is needed for the
application to  generating Somos-like sequences.

$\bullet$ {\tt FindSMpG(c,p,n,Maxr,d,a);} inputs a positive integer {\tt c}, a polynomial expression
{\tt p} in the discrete variable {\tt n}, positive integers  {\tt Maxr} and {\tt d}, and a symbol {\tt a}, and
outputs a non-linear recurrence, of order {\tt Maxr} and degree $\leq d$ satisfied by the sequence $a(n):=f(p(n))$, where
$f(n)$ is the sequence  satisfying $f(n)=cf(n-1)-f(n-2)$ and $f(0)=0, f(1)=1$. If it fails it returns
{\tt FAIL}, and the human has to try again with a larger {\tt d} and/or {\tt Marx}.
The first amazing example above was gotten from

{\tt FindSMpG(3,n**2,n,5,6,a);}

$\bullet$ {\tt FindSMeG(c,e,n,Maxr,d,a);} inputs a positive integer {\tt c}, a positive integer 
{\tt e},  positive integers  {\tt Maxr} and {\tt d} and a symbol {\tt a}, and
outputs a non-linear recurrence of order {\tt Maxr} and degree $\leq d$ satisfied by the sequence $f(e^n)$, where
$f(n)$ is the sequence  satisfying $f(n)=cf(n-1)-f(n-2)$ and $f(0)=0, f(1)=1$. If it fails it returns {\tt FAIL}.

The second amazing example above was gotten from:

{\tt FindSMeG(3,2,n,3,3,a);}

{\bf Webbbooks}  

For many amazing Somos-like miracles see the three webbooks

$\bullet$ {\tt http://www.math.rutgers.edu/\~{}zeilberg/tokhniot/oNesSomos1}

$\bullet$ {\tt http://www.math.rutgers.edu/\~{}zeilberg/tokhniot/oNesSomos3}

$\bullet$ {\tt http://www.math.rutgers.edu/\~{}zeilberg/tokhniot/oNesSomos4}

Enjoy!

Using {\tt NesSomos} you can find many more Somos-like {\it miracles} on your own,
impressing your friends and challenging your enemies.

{\bf Laurent Phenomenon}

Analogous things can be done to generate sequences that have the so-called Laurent phenomenon,
and even, {\it a priori}, the much stronger property of {\it positivity}.

{\bf Encore}

By looking at the output of  \hfill\break
{\tt http://www.math.rutgers.edu/\~{}zeilberg/tokhniot/oNesSomos1} \quad,
we see a pattern that can be summarized as follows.

{\bf Proposition}: Let $c$ be a positive integer (a formal symbol), and define a sequence  by
$$
a(0)=1 \quad , \quad a(1)=c \quad ,
$$
and for $n \geq 2$,
$$
a(n)=\frac{a(n-1)\, \left ( 4+(c-2)(c+2)a(n-1)^2 \, \right )}{2 \, + \, (c-2)(c+2)a(n-2)^2} \quad.
$$
Then for all $n \geq 1$, $a(n)$ are integers (polynomials in $c$ with integer coefficients).

{\bf Proof}: Consider the sequence defined by $f(0)=0, f(1)=1$ and for $n \geq 2$ by
$$
f(n)=cf(n-1)-f(n-2) \quad .
$$ 
Define $b(n):=f(2^n)$. It is routine to prove that 
$$
b(n)=\frac{b(n-1)(4+(c-2)(c+2)b(n-1)^2)}{2+(c-2)(c+2)b(n-2)^2} \quad ,
$$
and, of course,  $b(0)=1$ and $b(1)=c$. Hence, by
induction, $b(n)=a(n)$ for all $n \geq 1$. Since the $f(n)$ are obviously integers
(polynomials in $c$ with integer coefficients), so are $a(n)$. \halmos

{\bf References}

[K] Manuel Kauers,
{\it SumCracker: A package for manipulating symbolic sums and related objects},
Journal of Symbolic Computation {\bf 41} (2006),  1039-1057,
\hfill\break
{\tt http://www.risc.uni-linz.ac.at/people/mkauers/publications/kauers06h.pdf} \quad .

[KP] Manuel Kauers  and Peter Paule, {\it ``The Concrete Tetrahedron''}, Springer, 2011.

[M] J. L. Malouf, ``An integer sequence from a rational recursion'', Discr. Math. {\bf 110} (1992), 257-261.

[S] Neil Sloane, ``{\it The On-Line Encyclopedia of Integer Sequences}'', {\tt  http://oeis.org/} .

[Z1] Doron Zeilberger,
{\it An Enquiry Concerning Human (and Computer!) [Mathematical]
Understanding} in: C.S. Calude ,ed., ``Randomness and Complexity, from
Leibniz to Chaitin'', World Scientific, Singapore, Oct. 2007,
\hfill\break
{\tt http://www.math.rutgers.edu/\Tilde zeilberg/mamarim/mamarimhtml/enquiry.html} \quad .

[Z2] Doron Zeilberger,
{\it The C-finite Ansatz}, to appear in the Ramanujan Journal.
\hfill\break
{\tt http://www.math.rutgers.edu/\Tilde zeilberg/mamarim/mamarimhtml/cfinite.html} \quad .

\end